\documentclass{amsart}
\numberwithin{equation}{section}

\newtheorem{thm}{Theorem}[section]

\newtheorem{rem}[thm]{Remark}
\newtheorem{prob}[thm]{Problem}

\begin{document}

\begin{center}
\textbf{{\large {\ ON THE UNIQUENESS OF SOLUTIONS OF TWO INVERSE PROBLEMS FOR THE SUBDIFFUSION EQUATION}}}\\[0pt]
\medskip \textbf{Ravshan Ashurov$^{1}$ and Yusuf Fayziev$^{2}$}\\[0pt]
\medskip \textit{\ $^{1}$ Institute of Mathematics, Academy of Science of Uzbekistan}

\textit{ashurovr@gmail.com\\[0pt]
$^{2}$ National University of Uzbekistan}

\textit{fayziev.yusuf@mail.ru\\[0pt]
}
\end{center}

\textbf{Abstract}: Let $A$ be  an arbitrary positive selfadjoint
operator, defined in a
separable Hilbert space $H$.   The inverse
problems of determining the right-hand side of the equation and
the function $\varphi$ in the non-local boundary value problem $D_t^\rho u(t) + Au(t) =
f(t)$ ($0<\rho<1$, $0<t\leq T$), $u(\xi)=\alpha u(0)+\varphi$
($\alpha$ is a constant and  $0<\xi\leq T$), is considered. Operator $D_t$ on the left-hand
side of the equation expresses the Caputo derivative. For both inverse problems $u(\xi_1)=V$  is taken as the over-determination condition.  Existence and uniqueness theorems for
solutions of the problems under consideration are proved. The
influence of the constant $\alpha$ on the existence and uniqueness of a solution to problems is investigated. An interesting effect was discovered: when solving the forward problem, the uniqueness of the solution $u(t)$ was violated, while when solving the inverse problem for the same values of $\alpha$, the solution $u(t)$ became unique.

\textbf{Keywords}: Non-local problems, the
Caputo derivatives, subdiffusion equation, inverse problems.

\section{Introduction}

Let $A: H\rightarrow H$
be an arbitrary unbounded positive selfadjoint operator in a separable Hilbert space $H$ with the scalar product
$(\cdot, \cdot)$ and the norm $||\cdot||$. Let $A$ have a complete in $H$ system of orthonormal
eigenfunctions $\{v_k\}$ and a countable set of positive
eigenvalues $\lambda_k:$ $0<\lambda_1\leq\lambda_2 \cdot\cdot\cdot\rightarrow +\infty$. We will also assume that the sequence $\{\lambda_k\}$ has no finite limit points.
For a vector-valued functions (or simply functions)
$h: \mathbb{R}_+\rightarrow H$, we define the Caputo fractional derivative of order $0<\rho< 1$ as (see, e.g. \cite{Liz})
\[
D_t^\rho h(t)=\frac{1}{\Gamma
(1-\rho)}\int\limits_0^t\frac{h'(\xi)}{(t-\xi)^{\rho}} d\xi,
\quad t>0,
\]
provided the right-hand side exists. Here $\Gamma(\sigma)$ is
Euler's gamma function. Finally, let $C((a,b); H)$ stand
for a set of continuous functions $u(t)$  of $t\in (a,b)$ with
values in $H$.

The main object studied in this work is the following non-local boundary
value problem:
\begin{equation}\label{prob1}
\left\{
\begin{aligned}
&D_t^\rho u(t) + Au(t) = f(t),\quad 0< t \leq T;\\
&u(\xi_0) = \alpha u(0)+ \varphi, \quad 0 < \xi_0 \leq T,
\end{aligned}
\right.
\end{equation}
where $f(t) \in C((0,T]; H)$, $\varphi \in H$ and $\alpha$ is a
constant, $\xi_0 $ - a fixed  point. This problem is also called \emph{the forward
problem}.

In the case when $\xi_0=T$ and  parameter  $\alpha$ is equal to zero: $\alpha=0$,  this problem is called \emph{the backward problem} and it is well studied in the works \cite{Yama10} - \cite{Florida} and \cite{AA1}. And if $\alpha=0$ and $ \rho = 1 $, then we get a classical problem called the inverse heat
conduction problem with inverse time (\emph{retrospective inverse
problem}), which has been studied in detail by various specialists (see, e.g. Chapter 8.2 of \cite{Kab1} and literature therein).

It is well known that in most models described by differential (and pseudodifferential, see e.g., \cite{Umar}) equations, an initial condition is used to select a single solution. However, there are also processes where we have to use non-local conditions, for example, the integral over time intervals (see, e.g. \cite{Pao1} for reaction diffusion equations or \cite{Tuan} for fractional equations), or connection of solution values at different times, for example, at the initial time and at the final time (see, e.g. \cite{AshSob} - \cite{AshSob1}). It should be noted that non-local conditions model some details of natural phenomena more accurately, since they take into account additional information in the initial conditions.

The non-local boundary value problem (\ref{prob1}) for the classical diffusion equation, namely the following problem
\begin{equation}\label{probA}
	\left\{
	\begin{aligned}
		&u'(t) + Au(t) = f(t),\quad 0< t \leq T;\\
		&u(\xi_0) =u(0)+\varphi, \quad 0< \xi \leq T,
	\end{aligned}
	\right.
\end{equation}
has been extensively studied by many researchers
(see, e.g. A. O. Ashyralyev et al. \cite{AshSob} - \cite{AshSob1}). As shown in these papers, in contrast to the retrospective inverse problem,
problem (\ref{probA}) is coersively solvable in some spaces of
differentiable functions.

Let us return to the non-local problem (\ref{prob1}). The authors of this paper in their previous work \cite{AF2022} studied in detail the influence of parameter $\alpha\neq 0$ on the correctness of problem (\ref{prob1}). It turned out that the critical values of parameter $\alpha$ are in the interval $(0,1)$. In order to formulate the main result of work \cite{AF2022}, we recall the definition of the Mittag-Leffler function  $E_{\rho, \mu}(z)$ with two parameters (see, e.g. \cite{PSK}, Chapter 1):
\[
E_{\rho, \mu}(z)= \sum\limits_{n=0}^\infty \frac{z^n}{\Gamma(\rho
	n+\mu)},
\]
where $\mu$ is an arbitrary complex number. If parameter $\mu =1$, then we have the classical Mittag-Leffler function: $ E_{\rho}(z)= E_{\rho, 1}(z)$. Recall (see, e.g. \cite{AF2022}), $E_\rho(-t)$ decreases strictly monotonically as $t>0$ and, moreover, has the following estimate
\begin{equation}\label{ML}
	0<  E_\rho(-t) <1, \, t>0.
\end{equation}

In work \cite{AF2022} it is proved that if $\alpha \in (0,1)$ and $E_\rho(-\lambda_k t^\rho)\neq \alpha$ for all $k$, then the solution of problem (\ref{prob1}) exists and is unique. But it may turn out that for some eigenvalue $\lambda_{k_0}$ of operator $A$, with multiplicity $p_0$ (obviously, $p_0$ is a finite number), equality
\begin{equation}\label{critical}
    E_\rho(-\lambda_{k_0} t^\rho) = \alpha
\end{equation}
will hold. Then, as proved in \cite{AF2022}, in order for a solution to exist, it is necessary to require the following orthogonality conditions
\begin{equation}\label{Or1}
(\varphi,v_k)=0, \,\, (f(t), v_k)=0,\,\, \text{for all}\,\, t> 0,
\,\, k\in K_0;\,\, K_0=\{k_0, k_0+1,...., k_0+p_0-1\}.
\end{equation}
It should be noted that in this case there will be no uniqueness of the solution \cite{AF2022}.

The paper \cite{AF2022} also studies two inverse problems of determining the function $\varphi$ from the non-local condition (\ref{prob1}) and the source function $f$, i.e. the right-hand side of the equation in (\ref{prob1}) (in the latter case, it is assumed that $f$ does not depend on $t$). It is proved that if $\alpha \notin (0,1)$, then the solutions of both inverse problems exist and are unique. The main goal of this paper is to study these inverse problems for critical values of parameter $\alpha \in (0,1)$.

\begin{prob}\label{P1}
Let $\alpha \in (0,1)$. Find a pair $\{u(t), f \}$ of function $u(t)\in C([0,T]; H)$ and $f \in H$ with the properties
$D_t^\rho u(t), Au(t)\in C((0,T]; H)$ and satisfying the non-local problem (\ref{prob1}) (note, $f$ does not depend on $t$) and the over-determination condition
\begin{equation}\label{ODcon1}
u(\xi_1)=V, \quad 0 < \xi_1 <\xi_0,
\end{equation}
where $V$ is a given element of $H$.
\end{prob}
Note that if $\xi_1=\xi_0$, then the non-local condition in (\ref{prob1}) coincides with the Cauchy condition $u(0)=\varphi_1$ (note $\alpha \neq 0$). In this case, this inverse problem was studied in \cite{RuzTor}. If the reverse inequality $\xi_1 > \xi_0$ holds, then it will be shown that the solution may not be unique.

\begin{prob}\label{P2}
Let $\alpha \in (0,1)$. Find a pair $\{u(t), \varphi \}$ of function $u(t)\in C([0,T]; H)$ and $\varphi \in H$ with the properties
$D_t^\rho u(t), Au(t)\in C((0,T]; H)$ and satisfying the non-local problem (\ref{prob1}) and the over-determination condition
\begin{equation}\label{ODcon2}
u(\xi_2)=W, \quad 0 < \xi_2 \leq T, \,\, \xi_2 \neq \xi_0,
\end{equation}
where $W$ is a given element of $H$.
\end{prob}

If $\xi_2 = \xi_0$, then the non-local
condition $ u (\xi) = \alpha u (0) + \varphi $ coincides with the Cauchy condition $u(0)=\varphi_1$ (note $\alpha \neq 0$) and we have the backward problem,
considered in \cite{Yama10} - \cite{Florida}.

Everywhere below, for the vector - function $h(t)\in H$ (which may or may not depend on $t$) by the symbol $h_k(t)$ we will denote the Fourier coefficients with respect to the system of eigenfunctions $\{v_k\}$: $h_k(t)=(h(t), v_k)$.

\begin{thm}\label{thP1} Let  $\varphi, V \in D(A)$ and let the orthogonality conditions (\ref{Or1}) be satisfied. Then the inverse
Problem \ref{P1} has a unique solution
$\{u(t), f \}$ and this solution has the following form
\begin{equation}\label{fP1}
f =\sum\limits_{k\notin K_0} \bigg[
\frac{\alpha-E_{\rho}(-\lambda_k \xi_0^{\rho})}{ E_{\rho}(-\lambda_k
\xi_1^{\rho})\xi_0^\rho
E_{\rho,\rho+1}(-\lambda_k\xi_0^{\rho})+\xi_1^\rho
E_{\rho,\rho+1}(-\lambda_k \xi_1^{\rho})[\alpha-E_{\rho}(-\lambda_k
\xi_0^{\rho})]}\, V_k+
\end{equation}
$$+\frac{E_{\rho}(-\lambda_k \xi_1^{\rho})}{ E_{\rho}(-\lambda_k \xi_1^{\rho})\xi_0^\rho E_{\rho,\rho+1}
(-\lambda_k\xi_0^{\rho})+\xi_1^\rho E_{\rho,\rho+1}(-\lambda_k
\xi_1^{\rho})[\alpha-E_{\rho}(-\lambda_k \xi_0^{\rho})]}\,\varphi_k
\bigg] v_k,
$$
\begin{equation}\label{uP1}
u(t)= \sum\limits_{k\notin K_0} \left[\frac{E_{\rho}(-\lambda_k
t^{\rho})}{E_{\rho}(-\lambda_k \xi_0^{\rho})-\alpha}   \,
[{\varphi_k-f_k \xi_0^\rho E_{\rho,\rho+1}(-\lambda_k \xi_0^{\rho}) }]
+f_k t^\rho E_{\rho,\rho+1}(-\lambda_k t^{\rho})\right] v_k+.
\end{equation}
\[
+\sum\limits_{k\in K_0} \frac{E_\rho(-\lambda_k t^\rho)\,V_k}{E_\rho(-\lambda_k \xi_1^\rho)}\, v_k.
\]
\end{thm}
Note that, due to the orthogonality condition (\ref{Or1}), all Fourier coefficients $f_k$ vanish for $k\in K_0$. Obviously, $K_0$ can also be an empty set; in this case the sum $\sum_{k\notin K_0}$ is the same as $\sum_{k=1}^\infty$.

Let $ \tau $ be an arbitrary real number. In order to formulate a result on Problem \ref{P2}  we introduce the power
of operator $ A $, acting in $ H $ as
$$
A^\tau h= \sum\limits_{k=1}^\infty \lambda_k^\tau h_k v_k,
$$
where $h_k$ are the Fourier coefficients of $h\in H$. Obviously, the domain of this operator has the
form
$$
D(A^\tau)=\{h\in H:  \sum\limits_{k=1}^\infty \lambda_k^{2\tau}
|h_k|^2 < \infty\}.
$$
For elements of $D(A^\tau)$ we introduce the norm
\[
||h||^2_\tau=\sum\limits_{k=1}^\infty \lambda_k^{2\tau} |h_k|^2 =
||A^\tau h||^2,
\]
and together with this norm $D(A^\tau)$ turns into a Hilbert
space.

\begin{thm}\label{thP2} Let  $W \in D(A)$, $f \in C([0,T];D(A^\varepsilon))$ for some $\varepsilon\in (0,1)$
and let the orthogonality conditions (\ref{Or1}) be satisfied. Then the inverse
Problem \ref{P2} has a unique solution $\{u(t), \varphi\}$ and this solution has the form
\begin{equation}\label{fiP2}
\varphi= \sum\limits_{k\notin K_0} \left[\frac{E_{\rho}(-\lambda_k
\xi_0^{\rho})-\alpha}{E_{\rho}(-\lambda_k \xi_2^{\rho})} [W_k
-\omega_k(\xi_2)]+\omega_k(\xi_0)\right]  v_k,
\end{equation}
\begin{equation}\label{uP2}
u(t)= \sum\limits_{k\notin K_0} \left[\frac{\varphi_k-\omega_k
(\xi_0)}{E_{\rho}(-\lambda_k \xi_0^{\rho})-\alpha}   \,
E_{\rho}(-\lambda_k t^{\rho})+\omega_k(t)\right] v_k+
\end{equation}
\[
+\sum\limits_{k\in K_0}\frac{E_\rho(-\lambda_k t^\rho)\,W_k}{E_\rho(-\lambda_k \xi_2^\rho)}\, v_k,
\]
where
$$
\omega_k(t)=\int\limits_{0}^t\eta^{\rho-1}E_{\rho,\rho}(-\lambda_k
\eta^{\rho})f_k(t-\eta)d\eta.
$$
\end{thm}
Note that for $k\in K_0$ all Fourier coefficients $\varphi_k$ are equal to zero since the orthogonality condition (\ref{Or1}). When $K_0$ is an empty set, then the sum $\sum_{k\notin K_0}$ coincides with $\sum_{k=1}^\infty$.

\begin{rem}\label{unique}
 It should be specially noted that, as was proved in \cite{AF2022} and noted above, when equality (\ref{critical}) holds, the solution to the forward problem is not unique. But it turns out that both inverse problems have a unique solution even under condition (\ref{critical}).
\end{rem}

To the best of our knowledge,   the inverse problem of defining the function $ \varphi $ in the non-local condition was discussed  only  in the paper \cite{YulKad}. The authors considered this problem for the subdiffusion equation with the Caputo fractional derivative, the elliptic part of which is a two-variable differential expression with constant coefficients. On the other hand, it is not difficult to simulate a real process in which we will face just such an inverse problem. For example, in the temperature distribution process, the initial and final temperatures are not specified, and it is not required to find them, but information about the difference between the initial and final temperatures is sought.

As for the inverse problems of determining the source function $f$ with final time observation, it is well studied, both for classical partial differential equations and for equations of fractional order. Many theoretical studies have been published. Kabanikhin \cite{Kab1} and Prilepko, Orlovsky and Vasin \cite{prilepko} should be mentioned as classical monographs for integer-order equations. As for fractional differential equations, it is possible to construct theories parallel to the works of \cite{Kab1}, \cite{prilepko}, and work in this direction is ongoing. In this note, we will pay attention to only some of them, referring interested readers to a review paper \cite{Hand1}. Also note the works \cite{AF2022}, \cite{FurKir, LiYam, Sun}, where there is a review of recent work in this direction.

We note right away that no one has yet proposed a method for finding the right-hand side given in the abstract form $f(x,t)$. Known results deal with separated source
term $f(x, t) = q(t)p(x)$. The appropriate choice of the over-determination
depends on the choice whether the unknown is $q(t)$ or $p(x)$.

Quite a lot of papers are devoted to the case considered in this article, namely $q(t)\equiv 1$ and the unknown is $p(x)$. Subdiffusion equations whose elliptic part $A$ is an ordinary differential expression are considered, for example, in \cite{FurKir,KirMa, KirTor,TorTap}. The authors of \cite{LiYam, Sun},  studied the inverse problem for multi-term subdiffusion equations in which the elliptic part is either a Laplace operator or a second order operator. Article \cite{RuzTor} studied the inverse problem for the subdiffusion equation (\ref{prob1}) with the Cauchy condition. Recent articles \cite{AODif} - \cite{AOLob} are devoted to the inverse problem for the subdiffusion equation with Riemann-Liouville derivatives.

In  \cite{KirSal} non-self-adjoint differential operators (with non-local boundary conditions) were taken as A, and the solutions of the inverse problem were found in the form of a biorthogonal series.

In their previous work \cite{AF1}, the authors of this article considered an inverse problem for simultaneously determining the order of the Riemann-Liouville fractional derivative and the source function in the subdiffusion equations. Using the classical Fourier method, the authors proved the uniqueness and existence of a solution to this inverse problem.

It should be noted that in all of the listed works, the Cauchy
conditions in time are considered (an exception is work \cite{Saima}, where the integral condition is set with respect to the variable $t$). In the paper \cite{AF2022}, for the best of our knowledge, an inverse problem for subdiffusion equation with a non-local condition in time is considered for the first time.

The most difficult case to study is the case when the function $q(t)$ is unknown (see the survey paper \cite{Hand1} and \cite{Yama11} ). In inverse problems of this type, the condition $u(x_0,t)=u_0(t)$ is taken as an additional condition. The authors studied mainly the uniqueness of the solution of the inverse problem. In this regard, we note the recent papers \cite{ASh}, \cite{ASh2} where the inverse problem for determining the right-hand side of the form $q(t)$ was studied for the Schrodinger equation. Taking over-determination  conditions of a rather general form $Bu(\cdot, t)$, where $B: H\rightarrow R$ is a linear bounded functional, the authors proved both the existence and uniqueness of a solution to the inverse problem.

The papers \cite{AF2} - \cite{AF3} deal with the inverse problem of determining an order of the fractional derivative in the subdiffusion equation and in the wave equation, respectively.

\section{Inverse Problem \ref{P1}}

\textbf{2.1. Existence.} Assume that all the conditions of Theorem \ref{thP1} are satisfied, i.e. $\varphi, V \in D(A)$ and let the orthogonality conditions (\ref{Or1}) be satisfied.   Let us first prove the existence of a solution and that the solution has the form (\ref{fP1}) and (\ref{uP1}). The fact that  these series converge in the norm $H$ and in (\ref{uP1}) the summation and operators $D^\rho_t$ and $A$ can be interchanged was proved in the work of the authors \cite{AF2022}. Therefore, it suffices to show that the series (\ref{fP1}) and (\ref{uP1}) formally satisfy the equation and the initial condition (\ref{prob1}), and the over-determination condition (\ref{ODcon1}). In order to do this, we rewrite the series (\ref{fP1}) and (\ref{uP1}) in the form $f=\sum f_k v_k$ and $u(t)=\sum u_k(t) v_k$. Now, according to the Fourier method, it suffices to show that the unknown coefficients $f_k$ and $u_k(t)$ satisfy equation
\begin{equation}\label{Eq1}
    D_t^\rho u_k(t) + \lambda_k u_k(t) = f_k,
\end{equation}
the non-local condition
\begin{equation}\label{Nc1}
  u_k(\xi_0) = \alpha u_k(0) + \varphi_k,
    \end{equation}
and finally the over-determination condition
\begin{equation}\label{Ovc1}
  u_k(\xi_1) = V_k,
    \end{equation}
for all $k\geq 1$.

Let us show that $u_k(t)$ and $f_k$ satisfy equation (\ref{Eq1}). Let $k\notin K_0$. We have $u_k(t)=u_k^1(t)+ u_k^2(t)$, where
\[
 u_k^1(t)=\frac{E_{\rho}(-\lambda_k
t^{\rho})}{E_{\rho}(-\lambda_k \xi_0^{\rho})-\alpha}   \,
[{\varphi_k-f_k \xi_0^\rho E_{\rho,\rho+1}(-\lambda_k \xi_0^{\rho}) }]
    \]
    and
    \[
 u_k^2(t)=f_k t^\rho E_{\rho,\rho+1}(-\lambda_k t^{\rho}).
    \]
    It is known (see, e.g. \cite{Gor}, p. 174) that $u_k^1(t)$ is a solution to the homogeneous equation (\ref{Eq1}) with the initial condition $$
    u_k^1(0)=\frac{\varphi_k-f_k \xi_0^\rho E_{\rho,\rho+1}(-\lambda_k \xi_0^{\rho}) }{E_{\rho}(-\lambda_k \xi_0^{\rho})-\alpha}.
$$
It is also known (see ibid.) that the function
$$
\omega_k(t)=\int\limits_{0}^t\eta^{\rho-1}E_{\rho,\rho}(-\lambda_k
\eta^{\rho})f_k(t-\eta)d\eta
$$
from Theorem \ref{thP2} is a solution to equation (\ref{Eq1}) with the right-hand side $f_k(t)$ and with the initial condition $\omega_k(0)=0$. If in this formula $f_k(t)$ does not depend on $t$, then the integral can be rewritten in the form (see e.g. \cite{Gor}, formula (4.4.4))
\[
f_k\int\limits_{0}^t\eta^{\rho-1}E_{\rho,\rho}(-\lambda_k
\eta^{\rho})d\eta=f_k\,t^\rho E_{\rho,\rho+1}(-\lambda_k t^{\rho}).
\]
Therefore, the function $u_k^2(t)$ is a solution to the inhomogeneous equation (\ref{Eq1}) with the initial condition $ u_k^2(0)=0$.

Now suppose that $k\in K_0$. Then the function
\[
u_k(t)=\frac{E_\rho(-\lambda_k t^\rho)\,V_k}{E_\rho(-\lambda_k \xi_1^\rho)}
\]
is a solution of homogeneous equation (\ref{Eq1}) with the initial data
\[
u_k(0)=\frac{V_k}{E_\rho(-\lambda_k \xi_1^\rho)}.
\]
Thus, it is proved that the functions (\ref{fP1}) and (\ref{uP1}) really satisfy  equation (\ref{Eq1}).

It remains to verify the fulfillment of the non-local condition (\ref{Nc1}) and the over-determination condition (\ref{Ovc1}).

Let $k\notin K_0$. Since we have calculated $u_k(0)=u^1_k(0)+u^2_k(0)$, we can write
\[
\alpha u_k(0) +\varphi_k=\frac{\varphi_k E_\rho (-\lambda_k \xi_0^\rho)-\alpha f_k \xi_0^\rho E_{\rho, \rho+1}  (-\lambda_k \xi_0^\rho)}{E_\rho  (-\lambda_k \xi_0^\rho)-\alpha}.
\]
On the other hand, according to (\ref{uP1}),  $u_k(\xi_0)$ has exactly the same value:
\[
u_k(\xi_0) =\frac{\varphi_k E_\rho (-\lambda_k \xi_0^\rho)-\alpha f_k \xi_0^\rho E_{\rho, \rho+1}  (-\lambda_k \xi_0^\rho)}{E_\rho  (-\lambda_k \xi_0^\rho)-\alpha}.
\]

Let now $k\in K_0$. Then $\varphi_k=0$ (see (\ref{Or1})) and $E_\rho(-\lambda_k\xi_0^\rho)=\alpha$. Therefore
\[
\alpha u_k(0) +\varphi_k=\frac{\alpha V_k}{E_\rho(-\lambda_k \xi_1^\rho)},
\]
and
\[
u_k(\xi_0)=\frac{E_\rho(-\lambda_k \xi_0^\rho)\,V_k}{E_\rho(-\lambda_k \xi_1^\rho)}=\frac{\alpha V_k}{E_\rho(-\lambda_k \xi_1^\rho)}.
\]
Thus, the  Fourier coefficients of function $u(t)$, defined by formula (\ref{uP1}), satisfy the non-local condition (\ref{Nc1}) for all $k\geq 1$.

Let us check the fulfillment of the  over-determination condition (\ref{Ovc1}). Consider again the case $k\notin K_0$.
By virtue of  condition (\ref{Ovc1}) we obtain:
\[
 \frac{E_{\rho}(-\lambda_k \xi_1^{\rho})}{E_{\rho}(-\lambda_k \xi_0^{\rho})-\alpha}
 \, [{\varphi_k-f_k \xi_0^\rho E_{\rho,\rho+1}(-\lambda_k \xi_0^{\rho}) }] +f_k \xi_1^\rho E_{\rho,\rho+1}(-\lambda_k
 \xi_1^{\rho})=V_k.
\]
After simple calculations, we get
\[
f_k = \frac{\alpha-E_{\rho}(-\lambda_k \xi_0^{\rho})}{
E_{\rho}(-\lambda_k \xi_1^{\rho})\xi_0^\rho
E_{\rho,\rho+1}(-\lambda_k\xi_0^{\rho})+\xi_1^\rho
E_{\rho,\rho+1}(-\lambda_k \xi_1^{\rho})[\alpha-E_{\rho}(-\lambda_k
\xi_0^{\rho})]} V_k +
\]
$$
+\frac{E_{\rho}(-\lambda_k \tau^{\rho})}{ E_{\rho}(-\lambda_k
\tau^{\rho})\xi^\rho
E_{\rho,\rho+1}(-\lambda_k\xi^{\rho})+\tau^\rho
E_{\rho,\rho+1}(-\lambda_k \tau^{\rho})[\alpha-E_{\rho}(-\lambda_k
\xi^{\rho})]}\varphi_k,
$$
and this coincides with the Fourier coefficients of the function (\ref{fP1}).

If $k\in K_0$, then
\[
u_k(\xi_1)=\frac{E_\rho(-\lambda_k \xi_1^\rho)\,V_k}{E_\rho(-\lambda_k \xi_1^\rho)} = V_k.
\]

This completes the proof of the existence of a solution to Problem \ref{P1}.

\textbf{2.2. Uniqueness.} Let us proceed to the proof of the uniqueness of the solution of Problem \ref{P1}.

We proceed in the standard way: assuming the existence of two solutions, we obtain contradictions. Let $\{u_1(t),
f_1\}$ and $\{u_2(t), f_2\}$ be two solutions. It is required to prove $u(t)\equiv
u_1(t)-u_2(t)\equiv 0$ and $f\equiv f_1-f_2 = 0$. To determine $u(t)$ and $f$ we have the
problem:
\begin{equation}\label{eq1}
D_t^\rho u(t) + Au(t) = f, \quad  t>0;
\end{equation}
\begin{equation}\label{nl}
u(\xi_0) =\alpha u(0), \quad 0 < \xi_0 \leq T,
\end{equation}
\begin{equation}\label{od}
u(\xi_1) =0, \quad  0<\xi_1< \xi_0,
\end{equation}
where $\xi_0 $ and $\xi_1$ are the fixed  points.

Let  $u(t)$ be a solution to this problem and $u_k(t)=(u(t),
v_k)$. Then, by virtue of equation  (\ref{eq1}) and the
selfadjointness of operator $A$, problem (\ref{eq1})-(\ref{od}) becomes the following non-local problem with respect to $u_k(t)$:
\begin{equation}\label{prob2}
D_t^\rho u_k(t) +\lambda_k u_k(t)=f_k,\quad t>0; \quad u_k(\xi_0) =\alpha u_k(0), \quad u_k(\xi_1)=0.
\end{equation}
Note that if $k\in K_0$ then $f_k=0$.

Let first $k\notin K_0$. Suppose that $f_k$ is known and use the non-local condition to get (see, e.g. \cite{Gor}, p.174)
$$
u_k(t)=\frac{f_k \xi_0^{\rho}E_{\rho,\rho+1}(-\lambda_k \xi_0^{\rho})}{\alpha-E_{\rho}(-\lambda_k \xi_0^{\rho})}   \,
E_{\rho}(-\lambda_k t^{\rho})+f_k t^{\rho}E_{\rho,\rho+1}(-\lambda_k t^{\rho}).
$$
Now apply $u_k(\xi_1)=0$ to have
\begin{equation}\label{xitau}
f_k[ \xi_0^{\rho}E_{\rho,\rho+1}(-\lambda_k \xi_0^{\rho}) E_{\rho}(-\lambda_k \xi_1^{\rho})+ \xi_1^{\rho}E_{\rho,\rho+1}(-\lambda_k \xi_1^{\rho})(\alpha-E_{\rho}(-\lambda_k \xi_0^{\rho}))]=0.
\end{equation}
Let us show that for $\xi_1<\xi_0$ the square bracket is not equal to zero.
To do this, we introduce the notations: $a(t)=t^{\rho}E_{\rho,\rho+1}(-\lambda_k t^{\rho})>0$ and $b(t)=E_{\rho}(-\lambda_k t^{\rho})>0$. It is known (see, e.g. \cite{AF2022}) that the function $a(t)$ is increasing and the function $b(t)$ is decreasing. Now let us rewrite the square bracket as
\[
c(\xi_0, \xi_1)=a(\xi_0) b(\xi_1) - a(\xi_1) b(\xi_0) + \alpha b(\xi_0).
\]
Obviously, for $\xi_1<\xi_0$ this expression is strictly positive.
Therefore for all $k\notin K_0$ one has $f_k=0$ (see (\ref{xitau})).

It should be noted that if the inverse inequality $\xi_1 > \xi_0$ is satisfied, then the first term in the expression for $c(\xi_0, \xi_1)$ becomes less than the second one and, as a result, there is $\alpha\in (0,1)$ that turns $c(\xi_0, \xi_1)$ into zero. Therefore, in this case $f_k$ may not vanish, i.e., the uniqueness $f_k$ for these $\alpha$ and $k$ is violated.

Let us now consider the case $k\in K_0$. Denote $u_k(0)= b_k$. Then the unique solution to the
differential equation in (\ref{prob2}) with this initial condition has
the form $u_k(t)=b_k E_\rho (-\lambda_k t^\rho)$ (see, e.g. \cite{Gor}, p.174). Since $E_\rho (-\lambda_k \xi_0^\rho)=\alpha $ in the considering case, then the non-local condition is satisfied for an arbitrary $b_k$. But the over-determination condition $u_k(\xi_1)=0$ implies $b_k=0$ for $k\in K_0$.

Therefore, from the completeness of the system of eigenfunctions $ \{v_k \}
$, we finally obtain $f=0$ and $u(t) \equiv 0$, as required. The uniqueness and hence Theorem \ref{thP1} is completely proved.

\

\section{Inverse Problem \ref{P2}}

\textbf{3.1. Existence.} Suppose that $W \in D(A)$ and $f \in C([0,T];D(A^\varepsilon))$ for some $\varepsilon\in (0,1)$ and let the orthogonality conditions (\ref{Or1}) be satisfied.   Let us first show that series (\ref{fiP2}) and (\ref{uP2}) are indeed solutions to Problem \ref{P2}. The fact that $u(t)\in C([0,T]; H)$ and $\varphi \in H$ and have properties $D_t^\rho u(t), Au(t)\in C ((0,T]; H)$ was proved in our previous paper \cite{AF2022}.Therefore, it suffices to prove that (\ref{fiP2}) and (\ref{uP2}) together are a formal solution to Problem \ref{P2}. In turn, for this it suffices to show that the Fourier coefficients $\varphi_k$ and $u_k(t)$ of functions (\ref{fiP2}) and (\ref{uP2}) respectively,  satisfy equation (\ref{Eq1}), the non-local condition (\ref{Nc1}) and the over-determination condition
\begin{equation}\label{Ovc2}
u_k(\xi_2)=W_k.
    \end{equation}

It is not hard to verify that $u_k(t)$ is a solution of equation (\ref{Eq1}). Indeed, let first, $k\notin K_0$. We introduce the notation
\[
u_k^1(t)= \frac{\varphi_k-\omega_k
(\xi_0)}{E_{\rho}(-\lambda_k \xi_0^{\rho})-\alpha}   \,
E_{\rho}(-\lambda_k t^{\rho}).
\]
Then $u_k(t)=u_k^1(t)+ \omega_k(t)$. Here $u_k^1(t)$ is the solution of the homogeneous equation (\ref{Eq1}) with the initial condition
\[
u_k^1(0)= \frac{\varphi_k-\omega_k
(\xi_0)}{E_{\rho}(-\lambda_k \xi_0^{\rho})-\alpha},
\]
and $\omega_k(t)$ is the solution of equation (\ref{Eq1}) with zero initial condition (see, e.g. \cite{Gor}, p. 174).

If $k\in K_0$, then according to the orthogonality conditions $f_k=0$ and the function
\[
u_k(t)= \frac{E_\rho(-\lambda_k t^\rho)\,W_k}{E_\rho(-\lambda_k \xi_2^\rho)}
\]
is a solution of the homogeneous equation (\ref{Eq1}) with the initial condition
\[
u_k(0)= \frac{W_k}{E_\rho(-\lambda_k \xi_2^\rho)}.
\]

Thus we have shown that $u_k(t)$ is a solution of equation (\ref{Eq1}).

Let us check the non-local condition (\ref{Nc1}). Consider first the case $k\notin K_0$. We have
\[
\alpha u_k(0) +\varphi_k=\alpha\frac{\varphi_k-\omega_k
(\xi_0)}{E_{\rho}(-\lambda_k \xi_0^{\rho})-\alpha}+\varphi_k.
\]
On the other hand,
\[
u_k(\xi_0)=\frac{\varphi_k-\omega_k
(\xi_0)}{E_{\rho}(-\lambda_k \xi_0^{\rho})-\alpha}   \,
E_{\rho}(-\lambda_k \xi_0^{\rho})+\omega_k(\xi_0)=
\]
\[
=\frac{\varphi_k E_{\rho}(-\lambda_k \xi_0^{\rho}) -\alpha \, \omega_k
(\xi_0)}{E_{\rho}(-\lambda_k \xi_0^{\rho})-\alpha}=\frac{\varphi_k \big(E_{\rho}(-\lambda_k \xi_0^{\rho})-\alpha\big) +\varphi_k \alpha -\alpha \omega_k
(\xi_0)}{E_{\rho}(-\lambda_k \xi_0^{\rho})-\alpha}=
\]
\[
=\alpha\frac{\varphi_k-\omega_k
(\xi_0)}{E_{\rho}(-\lambda_k \xi_0^{\rho})-\alpha}+\varphi_k.
\]
Now consider the case $k\in K_0$. Note in this case $E_{\rho}(-\lambda_k \xi_0^{\rho})=\alpha$ and all Fourier coefficients $\varphi_k$ are equal to zero since the orthogonality condition (\ref{Or1}).
Therefore,
\[
\alpha u_k(0) +\varphi_k=\alpha\frac{W_k}{E_{\rho}(-\lambda_k \xi_2^{\rho})}=E_{\rho}(-\lambda_k \xi_0^{\rho})\frac{W_k}{E_{\rho}(-\lambda_k \xi_2^{\rho})}=u_k(\xi_0).
\]

Let us move on to checking the over-determination condition (\ref{Ovc2}). Let $k\notin K_0$. Then
$$
\frac{\varphi_k-\omega_k (\xi_0)}{E_{\rho}(-\lambda_k
\xi_0^{\rho})-\alpha} \, E_{\rho}(-\lambda_k
\xi_2^{\rho})+\omega_k(\xi_2) = W_k,
$$
or
$$
\varphi_k= \frac{E_{\rho}(-\lambda_k
\xi_0^{\rho})-\alpha}{E_{\rho}(-\lambda_k \xi_2^{\rho})} \,[W_k
-\omega_k(\xi_2)]+\omega_k(\xi_0),
$$
and this coincides with the Fourier coefficients of the function (\ref{fiP2}).

If $k\in K_0$, then
\[
u_k(\xi_2)=\frac{E_\rho(-\lambda_k \xi_2^\rho)\,W_k}{E_\rho(-\lambda_k \xi_2^\rho)} = W_k.
\]

This completes the proof of the existence of a solution to Problem \ref{P2}.

\textbf{3.2. Uniqueness.} Obviously, to prove the uniqueness of the solution to Problem \ref{P2}, it suffices to show that the solution $\{u(t), \varphi\}$  to the following inverse problem:
\[
D_t^\rho u(t) + Au(t) = 0, \quad \quad t>0;
\]
\[
u(\xi_0) =\alpha u(0)+ \varphi, \quad 0 < \xi_0 \leq T,
\]
\[
u(\xi_2) =0, \quad 0 < \xi_2 \leq T, \,\, \xi_2 \neq \xi_0,
\]
is identically zero: $u(t)\equiv 0$ and $\varphi=0$.

Let  $u(t)$ be a solution to this problem and let $u_k(t)=(u(t),
v_k)$. Then
\begin{equation}\label{Iprob2}
D_t^\rho u_k(t) +\lambda_k u_k(t)=0,\quad t>0; \quad u_k(\xi_0)
=\alpha u_k(0)+\varphi_k, \,\,  u_k(\xi_2) =0.
\end{equation}
Let $k\notin K_0$. Then it is not hard to verify that the  following function
$$
u_k(t)=\frac{E_{\rho}(-\lambda_k t^{\rho})}{E_{\rho}(-\lambda_k
\xi_0^{\rho})-\alpha}\,\, \varphi_k
$$
 is the only solution to the equation and non-local condition in (\ref{Iprob2}). The over-determination condition in (\ref{Iprob2}) implies
$$
u_k(\xi_2)=\frac{E_{\rho}(-\lambda_k
\xi_2^{\rho})}{E_{\rho}(-\lambda_k \xi_0^{\rho})-\alpha}\,\,
\varphi_k=0.
$$
Since $E_{\rho}(-\lambda_k \xi^{\rho})\neq\alpha$ and $E_{\rho}(-\lambda_k
\xi_2^{\rho})\neq 0$, then we have $\varphi_k=0$ and therefore $u_k(t)\equiv 0$ for all $k\notin K_0$.

Now consider the case  $k\in K_0$. Denote $u_k(0)= b_k$. Then the unique solution to the
differential equation in (\ref{Iprob2}) with this initial condition has
the form $u_k(t)=b_k E_\rho (-\lambda_k t^\rho)$ (see, e.g. \cite{Gor}, p.174). Since $E_\rho (-\lambda_k \xi_0^\rho)=\alpha $ and $\varphi_k=0$ in the considering case, then the non-local condition is satisfied for an arbitrary $b_k$. But the over-determination condition $u_k(\xi_2)=0$ implies $b_k=0$ and therefore $u_k(t)\equiv 0$  for $k\in K_0$.

Thus, from the completeness of the system of eigenfunctions $ \{v_k \}
$, we finally obtain $\varphi=0$ and $u(t) \equiv 0$, as required. The uniqueness and hence Theorem \ref{thP2} is completely proved.

\section{Conclusion}

In the previous paper of the authors \cite{AF2022} it is proved that for $\alpha \notin (0,1)$ the solutions of the forward and two inverse problems of determining $f$ and $\varphi$ exist and are unique. If $\alpha \in (0,1)$ and equality (\ref{critical}) holds for some $k\in K_0$, then to ensure the existence of the solution to the forward problem, it is necessary to require the orthogonality condition (\ref{Or1}). However, in this case the solution is not unique and it is determined up to the term
\[
\sum\limits_{k\in K_0} b_k E_\rho(-\lambda_k t^\rho)\, v_k,
\]
where $b_k$ are arbitrary numbers.

In this paper, we consider the above two inverse problems for critical values of parameter $\alpha \in (0,1)$. An interesting effect arises here: when solving the forward problem the uniqueness of solution $u(t)$ was violated, while when solving the inverse problem  for the same values of $\alpha$, solution $u(t)$ became unique. What is the matter here? It turns out, as follows from the main results of this paper, the over-determination condition
\[
u(\tau)=V
\]
can be rewritten in the form of two groups of conditions with respect to the Fourier coefficients
\[
u_k(\tau)=V_k, \,\,\, k\notin K_0,
\]
and
\[
u_k(\tau)=V_k, \,\,\, k\in K_0.
\]
With the help of the first group, the unique solutions of inverse problems are singled out, and since the coefficients $f_k$ and $\varphi_k$ are equal to zero for $k\in K_0$, the conditions from the second group are not used in this case. And the conditions from the second group ensure the uniqueness of the solution $u(t)$, namely, they determine uniquely the above arbitrary coefficients $b_k$.

\section{Acknowledgement}
The authors are grateful to Sh. A. Alimov for discussions of
these results.

The authors acknowledge financial support from the  Ministry of Innovative Development of the Republic of Uzbekistan, Grant No F-FA-2021-424.

\end{document}